\documentclass[3p]{elsarticle}

\usepackage{graphicx,color}
\usepackage{amsmath}
\usepackage{amssymb}
\makeatletter
\def\ps@pprintTitle{%
 \let\@oddhead\@empty
 \let\@evenhead\@empty
 \def\@oddfoot{}%
 \let\@evenfoot\@oddfoot}
\makeatother

\newcommand{\R}{\mathbb C}
\newtheorem{eg}{Example}[section]
\newtheorem{definition}[eg]{Definition}
\newtheorem{thm}[eg]{Theorem}
\newtheorem{lem}[eg]{Lemma}

\newtheorem{corollary}[eg]{Corollary}
\newenvironment{pf}{\noindent {\bf Proof:}}{\hspace*{\fill}$\square$}

\begin{document}

\begin{frontmatter}
\title{Zero Dynamics for Port-Hamiltonian Systems\footnote{The financial support  of the Oberwolfach Institute under the Research in Pairs Program, of the National Science and Engineering Research Council of Canada  Discovery Grant program, of the German Academic Exchange Service (DAAD) and of  the  Deutsche Forschungsgemeinschaft (DFG) for the research discussed in this article is gratefully acknowledged.  }} %
\author[BJ]{Birgit Jacob}
\ead{jacob@math.uni-wuppertal.de}
\address[BJ]{School of Mathematics and Natural Sciences, University of Wuppertal, Wuppertal, Germany}
 \author[KM]{Kirsten A.\ Morris}
 \ead{kmorris@uwaterloo.ca}
\address[KM]{Department of Applied Mathematics, University of Waterloo, Canada}
 \author[HZ,HZ2]{Hans Zwart}
 \ead{h.j.zwart@utwente.nl}
\address[HZ]{Department of Applied Mathematics, University of Twente, The Netherlands}
\address[HZ2]{Department of Mechanical Engineering, Eindhoven University of Technology, The Netherlands}

\date{\today}


\begin{abstract}                
The zero dynamics of infinite-dimensional systems can be difficult to characterize.  The zero dynamics of boundary control systems are particularly problematic. In this paper  the zero dynamics of  
 port-Hamiltonian systems are studied.  A complete characterization of the zero dynamics for
a port-Hamiltonian systems with invertible feedthrough as another port-Hamiltonian system on the same state space is given.  It is shown that the zero dynamics for any port-Hamiltonian system with commensurate wave speeds are well-defined, and are also   a port-Hamiltonian system. Examples include  wave equations with uniform wave speed on a network. A constructive procedure  for calculation of the zero dynamics, that can be used for very large system order,  is provided. 
\end{abstract}

\end{frontmatter}

\section{Introduction}

The zeros  of a  system are well-known to be important to controller design; see for instance, the textbooks \cite{DFT,Morristext}. For example,
the poles of a system controlled with a constant feedback gain move to the zeros
of the open-loop system as the gain increases. Furthermore, regulation is only possible if the
zeros of the system do not coincide with the poles of the signal to be tracked. Another
example is sensitivity reduction - arbitrary reduction of sensitivity is only possible if all the zeros are in the left half-plane. Right half-plane
zeros restrict the achievable performance; see for example, \cite{DFT}. The inverse of a system without right-hand-plane zeros can be approximated by a stable system, such systems are said to be minimum-phase and they are typically easier to control than non-minimum phase systems.

The zero dynamics are the dynamics of the system obtained by choosing the input $u$ so that the output $y$ is identically $0.$ This will  only be possible for initial conditions in some subspace of the original subspace. For linear systems with ordinary differential equation models and a minimal realization, the eigenvalues of the zero dynamics correspond to the zeros of the transfer function.  Zero dynamics are well understood for linear finite-dimensional systems, and have been extended to nonlinear finite-dimensional systems \cite{Isidori}.

However, many systems are modeled by delay or partial differential equations. This leads to an infinite-dimensional state space, and also an irrational transfer function. As for finite-dimensional systems, the zero dynamics are important.  
The notion of minimum-phase can be extended to infinite-dimensional systems; see in particular \cite{JMT} for a detailed study of conditions for second-order systems. Care needs to be taken since  a system can have no right-hand-plane zeros and still fail to be minimum-phase. The simplest such example  is a pure delay.
Results on adaptive control and on high-gain feedback
control of infinite-dimensional systems, see e.g.\
\cite{LogOwens,LogTown1,LogTown2,LogZwart,Nikitin}, require
the system to be minimum-phase. Moreover, the sensitivity of an
infinite-dimensional minimum-phase system can be reduced to an
arbitrarily small level   and  stabilizing controllers exist that achieve arbitrarily high gain
or phase margin \cite{FOT}.

Since the zeros are often not
accurately calculated by numerical approximations \cite{CM,Clark1997,GradMorris,Lindner}
it is useful to obtain an understanding of their behaviour in the
original infinite-dimensional context.
For infinite-dimensional systems with bounded control and observation, the  zero dynamics have been calculated, although they are not always well-posed \cite{MCSS,MorrisRebarberIJC,Zwartbook}. 

There are few results for zero dynamics for partial differential equations with boundary control and point observation.
In \cite{Byrnes,ByrnesGilliam}  the zero dynamics are found for a class of parabolic systems defined on an interval with collocated boundary control and observation. This was extended to the heat equation on an arbitrary region with collocated control and observation  in \cite{Reis}. 
In \cite{Kobayashi}  the invariant zeros for  a class of systems with analytic semigroup that includes boundary control/point sensing are defined and analysed.

 The zero dynamics of an important class of boundary control systems, port-Hamiltonian systems \cite{JZ,Gorrec:2005,Villegas_PhD} or systems of linear conservation laws \cite{Coron}, are  established in this paper. Such models are derived using a variational approach and many situations of interest, in particular waves and vibrations, can be described in a port-Hamiltonian framework.
A complete characterization of the zero dynamics for port-Hamiltonian systems with commensurate wave speeds is obtained. 
For any port-Hamiltonian system
 systems with invertible feedthrough, the  zero dynamics are another port-Hamiltonian system on the same state space. 
Port-Hamiltonian systems with commensurate wave speeds can be written as  as a coupling of scalar systems with the same wave speed. For these systems the zero dynamics are shown to be well-defined, and are in fact a  new port-Hamiltonian system. Preliminary versions of Theorem \ref{Tm:3} (for constant coefficients) and of Theorem \ref{T:4} (with an outline of the proof) appeared in \cite{IFAC2015}. 

A constructive procedure  for calculation of the zero dynamics based on linear algebra is provided.  This algorithm can be used on large networks, and does not use any approximation of the  system of partial differential equations.
The  results are illustrated with several examples.

\section{Port-Hamiltonian Systems}

Consider systems of the form  
\begin{eqnarray}
  \label{eq:2.1}
  \frac{\partial x}{\partial t}(\zeta,t) &=& P_1\frac{\partial}{\partial
\zeta}({\mathcal H}(\zeta)x(\zeta,t))
,\quad \zeta\in (0,1), t\ge 0 \\
x(\zeta,0)&=& x_0(\zeta),\quad \zeta\in (0,1)\\
  \label{eq:2.2}
  0 &=& W_{B,1}\begin{bmatrix} ({\mathcal H}x)(1,t)\\
    ({\mathcal H}x)(0,t)\end{bmatrix},\quad t\ge 0 \\
  \label{eq:2.3}
  u(t) &=& W_{B,2}\begin{bmatrix} ({\mathcal H}x)(1,t)\\
    ({\mathcal H}x)(0,t) \end{bmatrix},\quad t\ge 0\\
  \label{eq:2.4}
   y(t) &=& W_{C}\begin{bmatrix} ({\mathcal H}x)(1,t) \\
({\mathcal H}x)(0,t)\end{bmatrix},\quad t\ge 0,
\end{eqnarray}
where $P_1$ is an Hermitian  invertible $n\times n$-matrix, ${\mathcal H}(\zeta)$ is a positive $n\times n$-matrix for a.e.\  $\zeta\in (0,1)$ satisfying ${\cal H}, {\cal H}^{-1}\in L^\infty(0,1;\mathbb C^{n\times n})$, and $W_B := \left[\begin{smallmatrix}
    W_{B,1} \\ W_{B,2}\end{smallmatrix} \right]$ is a $n\times 2n$-matrix of rank $n$.
Such systems are said to be {\em port-Hamiltonian,} see \cite{Gorrec:2005,Villegas_PhD,JZ}, or systems of linear conservation laws \cite{Coron}. 
Here,  $x(\cdot,t)$ is the state of the system at time $t$, $u(t)$ represents the input of the system at time $t$ and $y(t)$ the output of the system at time $t$.

A different representation of port-Hamiltonian systems, the diagonalized form, will be used.
The matrices $P_1 {\mathcal H}(\zeta)$ possess the same eigenvalues counted according to their multiplicity as the matrix $ {\mathcal H}^{1/2}(\zeta) P_1 {\mathcal H}^{1/2} (\zeta)$, and as ${\mathcal H}^{1/2}(\zeta) P_1 {\mathcal H}^{1/2}(\zeta)$ is diagonalizable  the  matrix  $P_1 {\mathcal H}(\zeta)$ is diagonalizable as well. Moreover, by our assumptions, zero is not an eigenvalue  of  $P_1 {\mathcal H}(\zeta)$ and all eigenvalues are real, that
is, there exists an invertible matrix $S(\zeta)$ such that
\[ 
  P_1{\mathcal H}(\zeta) =S^{-1} (\zeta)\underbrace{{\rm diag} (p_1(\zeta), \cdots, p_k(\zeta), n_1(\zeta), \cdots , n_l(\zeta))}_{=:\Delta(\zeta)} S(\zeta).
\]
Here $p_1(\zeta), \cdots,p_k(\zeta)>0$ and $n_1(\zeta),\cdots, n_l(\zeta)<0$. In the remainder of this article  it is assumed that $S$ and $\Delta$ are continuously differentiable on $(0,1)$.
Introducing the new state vector 
\[ z(\zeta,t)=  \begin{bmatrix}  z_+(\zeta,t) \\ z_-(\zeta,t)\end{bmatrix}  =S(\zeta) x(\zeta,t), \qquad \zeta \in [0,1],\]
with $z_+(\zeta,t) \in \mathbb C^k$ and $z_-(\zeta,t) \in \mathbb C^l$,
and writing
\[ \Delta(\zeta) =\begin{bmatrix}    \Lambda(\zeta) & 0\\0 & \Theta (\zeta)  \end{bmatrix},\] 
where $\Lambda(\zeta)$ is a positive definite $k\times k$-matrix and $\Theta(\zeta)$ is a negative definite $l\times l$-matrix,
the system \eqref{eq:2.1}--\eqref{eq:2.4} can be equivalently written as
\begin{eqnarray}
\label{eq:2.1b}
  \frac{\partial }{\partial t} z(\zeta,t)   &=& \frac{\partial}{\partial
\zeta} \left(\Delta(\zeta)  z(\zeta,t)\right) +S(\zeta)\frac{S^{-1}(\zeta)}{d\zeta}\Delta(\zeta)z(\zeta,t),\\
z(\zeta,0)&=& z_0(\zeta),\quad \zeta\in (0,1)\\
    \label{eq:2.2b}
\begin{bmatrix}
0 \\ u (t) \end{bmatrix} &=&
\underbrace{\begin{bmatrix} K_{0+} & K_{0-} \\ K_{u+}  & K_{u-}  \end{bmatrix} }_{ K} \begin{bmatrix}  \Lambda(1) z_+(1,t)  \\ \Theta(0) z_-(0,t) \end{bmatrix} +
\underbrace{\begin{bmatrix} L_{0+} &  L_{0-}  \\ L_{u+}  &  L_{u-}  \end{bmatrix}  }_{ L} \begin{bmatrix}  \Lambda(0) z_+(0,t)  \\ \Theta(1) z_-(1,t) \end{bmatrix} , \\
  \label{eq:2.4b}
 y(t) &=& \underbrace{\begin{bmatrix}K_{y+} & K_{y-} \end{bmatrix}}_{K_y} \begin{bmatrix}  \Lambda(1) z_+(1,t)  \\ \Theta(0) z_-(0,t) \end{bmatrix} + \underbrace{ \begin{bmatrix}L_{y+} & L_{y-} \end{bmatrix}}_{L_y} \begin{bmatrix}  \Lambda(0) z_+(0,t)  \\ \Theta(1) z_-(1,t) \end{bmatrix},
\end{eqnarray}
where $t\ge 0$ and $\zeta\in(0,1)$. 

Next, consider well-posedness of the system \eqref{eq:2.1b}--\eqref{eq:2.4b}, or equivalently of system  \eqref{eq:2.1}--\eqref{eq:2.4}. Well-posedness means that for every initial condition $z_0 \in
L^2(0,1;{\mathbb C}^n)$ and every input $u\in
L^2_{\mathrm{loc}}(0,\infty;{\mathbb C}^p)$ the mild solution
$z$ of the
system (\ref{eq:2.1b})--(\ref{eq:2.2b}) is well-defined in the state
space $X:= L^2(0,1;{\mathbb C}^{n})$ and the output (\ref{eq:2.4b}) is well-defined in
$L^2_{\mathrm{loc}}(0,\infty;{\mathbb C}^m)$. See  \cite{JZ} for the precise definition and further results on well-posedness of port-Hamiltonian systems.
To characterize well-posedness, define the matrices
$$K=\begin{bmatrix} K_0\\K_u\end{bmatrix} = \begin{bmatrix} K_{0+} & K_{0-} \\ K_{u+}  & K_{u-}  \end{bmatrix}, \quad L=\begin{bmatrix} L_0\\L_u\end{bmatrix} = \begin{bmatrix} L_{0+} &  L_{0-}  \\ L_{u+}  &  L_{u-}  \end{bmatrix}.
$$
\begin{thm}
\label{thm_JZ} \cite{ZwartESAIM2010},
\cite[Thm.\ 13.2.2 and 13.3.1]{JZ}. The following are equivalent
\begin{enumerate}
\item The system \eqref{eq:2.1b}--\eqref{eq:2.4b} is well-posed on  $L^2(0,1;\mathbb C^{n})$;
\item For every initial condition $z_0 \in L^2(0,1;\mathbb C^{n})$, the partial differential equation (\ref{eq:2.1b})--(\ref{eq:2.2b}) with $u=0$ possesses a unique mild solution on the state space $L^2(0,1;\mathbb C^{n})$. Furthermore, this solution depends continuously on the initial condition;
\item The matrix $K$ is invertible. 
\end{enumerate}
\end{thm}

Thus, well-posedness of a port-Hamiltonian system is equivalent to well-posedness of a homogeneous partial differential equation; the boundedness of the input/state and state/output maps does not need to be checked separately. For the remainder of this paper it is assumed that  $K$ is invertible so that the control system is well-posed. The corresponding generator $A$ of the $C_0$-semigroup of the homogeneous system is given by \cite{JZ}
\begin{align*}
 Af=& -(\Delta f)' +S(S^{-1})'\Delta f,\\
 D(A)=&\left\{\Delta f\in H^1(0,1;\mathbb C^n)\mid \begin{bmatrix} 0 \\ 0 \end{bmatrix}  = K  \begin{bmatrix} \Lambda(1)f_+(1)  \\ \Theta(0)f_-(0)   \end{bmatrix}  + L
\begin{bmatrix} \Lambda(0)f_+(0)  \\ \Theta(1)f_-(1)   \end{bmatrix}  \right\}.
\end{align*}
The resolvent operator of $A$ is compact, and thus the spectrum of $A$ contains only eigenvalues.
For port-Hamiltonian systems, well-posedness implies that the system
(\ref{eq:2.1b})--(\ref{eq:2.4b}) is also regular, i.e, the transfer function $G(s)$ possesses a limit over the real line, see
\cite{ZwartESAIM2010} or \cite[Section 13.3]{JZ}. Writing 
\begin{equation} 
\label{feedthrough} 
  K_y K^{-1}= \begin{bmatrix} * & E \end{bmatrix} 
\end{equation}
with $E\in \mathbb C^{m\times p}$,  this limit of $G(s)$ over the real axis is $E$, see \cite[Theorem 13.3.1]{JZ}.

Now consider zero dynamics for  port-Hamiltonian systems. %
\begin{definition}
Consider the system (\ref{eq:2.1b})--(\ref{eq:2.4b}) on the state space $X=L^2(0,1;{\mathbb C}^n)$. 
The {\em zero dynamics} of  (\ref{eq:2.1b})--(\ref{eq:2.4b}) are the pairs $(z_0,u)\in X\times L^2_{\mathrm{loc}}(0,\infty;{\mathbb C^p})$
for which the mild solution of (\ref{eq:2.1b})--(\ref{eq:2.4b}) satisfies $ y = 0$.
The largest output nulling subspace is 
  \begin{eqnarray*}
    \nonumber
    V^* &=& \{ z_0 \in X \mid \mbox{ there exists a function } u \in L^2_{\mathrm{loc}}(0,\infty;{\mathbb C^p}): \\
  \label{eq:15}
   &&  \quad \mbox{the mild solution of (\ref{eq:2.1b})--(\ref{eq:2.4b}) satisfies } y = 0 \}.
  \end{eqnarray*}
\end{definition}
Thus, $V^*$ is the space of initial conditions for which there exists a control $u$ that ``zeroes'' the output.

Setting $y=0$ in (\ref{eq:2.4b})  reveals that the zero dynamics are described by 
\begin{eqnarray}
\label{eq:2.1z}
  \frac{\partial }{\partial t} z(\zeta,t)   &=& \frac{\partial}{\partial
\zeta} \left(\Delta(\zeta)  z(\zeta,t)\right) +S(\zeta)\frac{S^{-1}(\zeta)}{d\zeta}\Delta(\zeta)z(\zeta,t),\\
z(\zeta,0)&=& z_0(\zeta),\quad \zeta\in (0,1)\\
    \label{eq:2.2z}
0 &=&
\begin{bmatrix} K_{0}  \\ K_{y}  \end{bmatrix}   \begin{bmatrix}  \Lambda(1) z_+(1,t)  \\ \Theta(0) z_-(0,t) \end{bmatrix} +
\begin{bmatrix} L_{0}  \\ L_{y}  \end{bmatrix}   \begin{bmatrix}  \Lambda(0) z_+(0,t)  \\ \Theta(1) z_-(1,t) \end{bmatrix}, \\
  \label{eq:2.4z}
 u(t) &=& K_u \begin{bmatrix}  \Lambda(1) z_+(1,t)  \\ \Theta(0) z_-(0,t) \end{bmatrix}  + L_u \begin{bmatrix}  \Lambda(0) z_+(0,t)  \\ \Theta(1) z_-(1,t) \end{bmatrix},
\end{eqnarray}
where $t\ge 0$ and $\zeta\in(0,1)$. Note that system  \eqref{eq:2.1z}--\eqref{eq:2.4z} is still in the format of a port-Hamiltonian system, but even regarding  (\ref{eq:2.4z}) as the (new) output, it needs not to be a well-posed port-Hamiltonian system since the new ``$K$-matrix'', $\left[\begin{smallmatrix} K_{0}  \\ K_{y}  \end{smallmatrix} \right]$ can have rank less than $n$.  
The zero dynamics are a well-posed dynamical system if 
the system \eqref{eq:2.1z}--\eqref{eq:2.4z} with state-space $V^*$, no input and output $u$ is well-posed. 

The eigenvalues of the zero dynamics of the system are closely related to the invariant and transmission zeros of the system. For simplicity only the single-input single-output case is considered ($p=m=1$).
\begin{definition}\cite{Reis, CM}
A complex number $\lambda\in\mathbb C$ is  an {\em  invariant zero} of  the system (\ref{eq:2.1b})--(\ref{eq:2.4b}) on the state space $X=L^2(0,1;{\mathbb C}^n)$, if there exist  $z\in H^1(0,1;\mathbb C^n)$ and $u\in \mathbb C$ such that
\begin{eqnarray*}
  \lambda  z(\zeta)   &=& \frac{\partial}{\partial
\zeta} \left(\Delta(\zeta)  z(\zeta)\right) +S(\zeta)\frac{S^{-1}(\zeta)}{d\zeta}\Delta(\zeta)z(\zeta),\\
0 &=&
\begin{bmatrix} K_{0}  \\ K_{y}  \end{bmatrix}   \begin{bmatrix}  \Lambda(1) z_+(1)  \\ \Theta(0) z_-(0) \end{bmatrix} +
\begin{bmatrix} L_{0}  \\ L_{y}  \end{bmatrix}   \begin{bmatrix}  \Lambda(0) z_+(0)  \\ \Theta(1) z_-(1) \end{bmatrix}, \\
 u &=& K_u \begin{bmatrix}  \Lambda(1) z_+(1)  \\ \Theta(0) z_-(0,t) \end{bmatrix}  + L_u \begin{bmatrix}  \Lambda(0) z_+(0,)  \\ \Theta(1) z_-(1,t) \end{bmatrix},
\end{eqnarray*}
\end{definition}

\begin{definition}
A complex number $s\in\mathbb C$ is  a {\em transmission zero} of  the system (\ref{eq:2.1b})--(\ref{eq:2.4b})  if the transfer function satisfies $G(s)=0$.
\end{definition}

If $\lambda \in \rho(A)$, where $\rho(A)$ denotes the resolvent set of $A$, then $\lambda$ is an invariant zero if and only if $\lambda$ is a transmission zero \cite[Theorem 12.2.1]{JZ}. Moreover, if the zero dynamics is well-posed, then the spectrum of the corresponding generator equals the set of invariant zeros of the  system (\ref{eq:2.1b})--(\ref{eq:2.4b}).

If the feedthrough operator of the original system is invertible, then  the zero dynamics system is well-posed on the entire state space, and is also a port-Hamiltonian system.  
\begin{thm}
\label{Tm:3}
Assume that the system has the same number of inputs as outputs. Then the zero dynamics are well-posed on the entire state space if and only if the feedthrough operator $E$ of the original system is invertible.
\end{thm}
{\bf Proof:} 
This was proven in \cite{IFAC2015} in the case of a constant coefficient matrix $\mathcal H$. The proof presented here 
is more complete, and 
includes the generalization to variable coefficients.
The  feedthrough operator $E$ of the original system is given by  $[* ~ E] = K_y  K^{-1}$ (see \eqref{feedthrough}). 
It will first be shown 
that  invertibility of $E$ is equivalent to invertibility of the ``$K$-matrix'' of equation (\ref{eq:2.2z}):
 $$
  \tilde K:= \begin{bmatrix} K_{0}  \\ K_{y} \end{bmatrix}.
$$
 
If $E$ is singular, then there is $u \neq 0$ in the kernel of $E$,  and
\[ K_y  K^{-1}\begin{bmatrix}
    0 \\ u \end{bmatrix} =0.
\]
Combining this with the fact that $K_0 K^{-1} = \begin{bmatrix} I &
  0 \end{bmatrix}$, 
\[ \tilde KK^{-1} \begin{bmatrix}
    0 \\ u \end{bmatrix} =\begin{bmatrix} K_0 \\ K_y  \end{bmatrix}  K^{-1} \begin{bmatrix}
    0 \\ u \end{bmatrix} =0.
\]
Thus $\tilde{K}$ is singular. Assume next that $\tilde{K}$ is singular. Thus there exists non-zero 
$\left[\begin{smallmatrix} x_1 \\ x_2 \end{smallmatrix}\right]$ such that
\begin{equation}
\label{eq:1HZ}
  \begin{bmatrix} K_0 \\ K_y  \end{bmatrix}  \begin{bmatrix}
    x_1 \\ x_2 \end{bmatrix} = \begin{bmatrix}
    0 \\ 0 \end{bmatrix} .
\end{equation}
This implies that 
\[
   K \begin{bmatrix}
    x_1 \\ x_2 \end{bmatrix} = \begin{bmatrix} K_0 \\ K_u \end{bmatrix}  \begin{bmatrix}
    x_1 \\ x_2 \end{bmatrix} = \begin{bmatrix}
    0 \\ z \end{bmatrix} ,
\]
where $z\neq 0$, since $K$ is invertible. Thus
\[
 Ez =  K_y  K^{-1}\begin{bmatrix}
    0 \\ z \end{bmatrix}= K_y \begin{bmatrix}
    x_1 \\ x_2 \end{bmatrix} =0
\]
and thus $E$ is not invertible.

Assume now that $E$ is invertible, then by the above equivalence with the invertibility of $\tilde{K}$ and Theorem \ref{thm_JZ} for every initial condition there exists a solution of (\ref{eq:2.1z})--(\ref{eq:2.2z}). Since $z$ is now determined,  $u$ is determined by  (\ref{eq:2.4z}). Now it is straightforward to see that the functions $z$ and $u$  satisfy (\ref{eq:2.1b})--(\ref{eq:2.2b}) and the corresponding output $y$ satisfies $y=0$.

If for every $z_0\in L^2(0,1;{\mathbb C}^n)$ there exists a solution of (\ref{eq:2.1z})--(\ref{eq:2.4z}), then the functions $z$ and $u$  satisfy (\ref{eq:2.1b})--(\ref{eq:2.2b}). Since $K$ is invertible, the solution depends continuously on the initial condition. By construction, $z$ is the solution of the homogeneous equation (\ref{eq:2.1z})--(\ref{eq:2.2z}), and  Theorem \ref{Tm:3} implies the invertibilty of $\tilde{K}$.
\hfill $\square$
\medskip

Example 1 in \cite{IFAC2015} illustrates calculation of the zero dynamics in the case where $E$ is invertible. 

It is very common though for the feedthrough to be non-invertible. This more challenging  situation is considered in the next two sections.

\section{Commensurate constant wave speed}
\label{sec:3}

In this section, the following class of port-Hamiltonian systems is considered:
\begin{eqnarray}
\label{eq:2.1d}
  \frac{\partial }{\partial t}  z(\zeta,t)   &=&  -\lambda_0\frac{\partial}{\partial \zeta}  z(\zeta,t),\\
z(\zeta,0)&=& z_0(\zeta),\quad \zeta\in (0,1)\\
    \label{eq:2.2d}
\begin{bmatrix}
0 \\ u (t) \end{bmatrix} &=& -\lambda_0 K   z(0,t) -\lambda_0 L  z(1,t) , \\
  \label{eq:2.4d}
 y(t) &=& -\lambda_0 K_{y}  z(0,t)  -\lambda_0 L_{y}  z(1,t).
\end{eqnarray}
If ${\mathcal H}$ is constant, then  (\ref{eq:2.1b})--(\ref{eq:2.4b}) is of the form  (\ref{eq:2.1d})--(\ref{eq:2.4d}) with $-\lambda_0$ replaced by a diagonal (constant) and invertible matrix $\Delta$.
On the diagonal of the matrix  $\Delta$ are the possible different wave speeds of the system. If the ratio of any pair of diagonal entries of $\Delta$ is rational, then the system \eqref{eq:2.1b}--\eqref{eq:2.4b} can be equivalently written in form  \eqref{eq:2.1d}--\eqref{eq:2.4d}  by dividing the intervals to adjust the propagation periods. This is a standard procedure and is illustrated in Example \ref{E:3.1}. The following simple reflection makes positive wave speeds into negative wave speed, while keeping the same absolute speed
\[
  \tilde{z}_k(\zeta,t) :=  z_k(1-\zeta,t).
\]
It is good to remark that the system \eqref{eq:2.1d}--\eqref{eq:2.4d} will in general have larger matrices than the original system (\ref{eq:2.1b})--(\ref{eq:2.4b}). However, for simplicity, still denote the size  by $n$. 
\begin{eg}
\label{E:3.1}
Consider the following system with commensurable wave speeds
\begin{eqnarray*}
 \frac{ \partial z_1}{\partial t} &=& - \frac{ \partial z_1}{\partial \zeta }, \\
  \frac{ \partial z_2}{\partial t} &=&- \frac{1}{2} \frac{ \partial z_2}{\partial \zeta }, 
\end{eqnarray*}
with $\zeta\in [0,1]$, $t\ge 0$ and
\begin{eqnarray*}
\begin{bmatrix} 0 \\ u (t) \end{bmatrix} &=& \begin{bmatrix} 1 & 1  \\ 0 & 1   \end{bmatrix} z(0,t) + \begin{bmatrix} 1 & 0  \\ 0 &  0 \end{bmatrix} z(1,t) \\
y(t) & = & \begin{bmatrix} 0 & 0  \end{bmatrix} z(0,t) + \begin{bmatrix} 1 & 1 \end{bmatrix} z(1,t) . 
\end{eqnarray*}
This system has not a uniform wave speed, but can be written equivalently as a system with one wave speed. To reach this goal, split  the second equation in two and obtain the following equivalent system
\begin{eqnarray*}
 \frac{ \partial z_1}{\partial t} &=& - \frac{ \partial z_1}{\partial \zeta }, \\
  \frac{ \partial z_{2a}}{\partial t} &=& - \frac{ \partial z_{2a}}{\partial \zeta }, \\
   \frac{ \partial z_{2b}}{\partial t} &=&  -\frac{ \partial z_{2b}}{\partial \zeta }, 
\end{eqnarray*}
with $\zeta\in [0,1]$, $t\ge 0$, $z_{2b}(\zeta,t) = z_2(\zeta/2,t)$ and $z_{2a}(\zeta,t) = z_2((1+\zeta)/2,t)$ and
\begin{eqnarray*}
\begin{bmatrix} 0 \\ 0 \\ u (t) \end{bmatrix} &=& \begin{bmatrix} 1 & 0 &1  \\ 0 & 1& 0 \\ 0&0 &1  \end{bmatrix} z(0,t) + \begin{bmatrix} 1 & 0 &0 \\ 0 &  0 &-1 \\ 0 &0&0\end{bmatrix} z(1,t) \\
y(t) & = & \begin{bmatrix} 0 & 0 &0 \end{bmatrix} z(0,t) + \begin{bmatrix} 1 & 1 &0\end{bmatrix} z(1,t) . 
\end{eqnarray*}
\end{eg}
This transformation also works  if ${\cal H}(\zeta)$ is diagonal a.e.\ $\zeta\in(0,1)$ and the ratio of the numbers $\tau_i:=\int_0^1 \frac{1}{{\cal H}(\zeta)_{ii}}d\zeta$ are pairwise rational  \cite{suimai}. 

It  is now shown that  the zero dynamics can be well-defined through the input and output equations.

It is well-known that the solution of \eqref{eq:2.1d} is given by $ z(\zeta,t)= f(  1- \zeta + \lambda_0 t )$
for $t\ge 0$ and some function $f$. Using this fact,  we write the system \eqref{eq:2.1d}-\eqref{eq:2.4d}
equivalently as
\begin{eqnarray}
f(t) &=& z_0(1-t),\qquad t\in [0,1],\label{eq:17a}\\
 \begin{bmatrix} 0\\ u(t) \end{bmatrix} &=& \label{eq:17}
   -\lambda_0 K   f(  1+\lambda_0 t )  -\lambda_0 L   f(  \lambda_0  t ),\quad t\ge 0,\\
  \label{eq:22}
  y(t) &=&  -\lambda_0K_y   f(1+\lambda_0 t)  -\lambda_0 L_y f( \lambda_0 t ),\quad t\ge 0.
\end{eqnarray}
Since the system is well-posed,  the matrix $K$ is invertible (Theorem \ref{thm_JZ}). 
Thus, equivalently
\begin{eqnarray}
f(t) &=& z_0(1-t),\qquad t\in [0,1],\label{eq:18a}\\
  \label{eq:18}
   f( 1+ \lambda_0 t)  &=& - K^{-1} L  f( \lambda_0 t) -\lambda_0^{-1}K^{-1} \begin{bmatrix} 0 \\u(t)  \end{bmatrix},\quad t\ge 0, \\
   y(t) &=& (\lambda_0K_y K^{-1} L-\lambda_0L_y) f(\lambda_0 t) + K_y K^{-1} \begin{bmatrix} 0 \\u(t)  \end{bmatrix},\quad t\ge 0. \label{eq:18b}
\end{eqnarray}
Defining
\begin{eqnarray}
  A_d&=&-  K^{-1} L,\, \quad B_d=  - \lambda_0^{-1} K^{-1} \begin{bmatrix} 0 \\ I \end{bmatrix},  \nonumber \\
 C_d&=&-\lambda_0 K_yA_d -\lambda_0 L_y, \quad D_d=-\lambda_0 K_yB_d,  \label{disc}
\end{eqnarray}
equation \eqref{eq:18}--\eqref{eq:18b} can be written as
\begin{eqnarray*}
  f(1+\lambda_0 t) &=& A_d f (\lambda_0 t) + B_d u(t),\\
  y(t) &=& C_d f(\lambda_0 t) +D_d u(t).
  \end{eqnarray*}
Define for $n \in {\mathbb N}$ the functions  $z_d (n) \in L^2 (0, 1;\R^n)$, $u_d (n)\in L^2 (0, 1;\R^p), $ and $y_d (n ) \in  L^2 (0, 1;\R^m)$ by $z_d(0)(\zeta):=z_0(1-\zeta)$, $z_d (n) (\zeta ) = f(n+ \zeta  ) $ for $n\ge 0$ and 
$$ u_d (n) (\zeta ) = u(\frac{n+\zeta}{\lambda_0}), \quad y_d (n) (\zeta ) = y(\frac{n+\zeta}{\lambda_0} ),\qquad n\in\mathbb N. $$
Thus equations \eqref{eq:2.1d}--\eqref{eq:2.4d} can be equivalently rewritten as
\begin{eqnarray}
\label{eq:19}
  z_d(n+1) (\zeta) &=& A_d z_d(n) (\zeta)  + B_d u_d(n) (\zeta) \\
 (z_d(0))(\zeta) &=& z_0(1-\zeta)\\
\label{eq:20}
  y_d(n) (\zeta) &=& C_d z_d (n) (\zeta) + D_d u_d(n) (\zeta)
\end{eqnarray}
This representation is very useful,  not only  for the zero dynamic, but also for other  properties like stability.

%
\begin{thm}
\label{thm32}\cite[Corollary 3.7]{kloess}
The system   \eqref{eq:2.1d}--\eqref{eq:2.4d} is exponentially stable if and only if  the spectral radius of $A_d$ satisfies $r (A_d) <1$ or equivalently if $\sigma_{max} (A_ d )<1$.
\end{thm}

Further sufficient conditions for exponential stability can be found in \cite{Coron, Engel, JZ}. In particular, exponential stability is implied by 
the  condition  $K K^* - L L^* >0$,  \cite[Thm. 3.2]{Coron} and \cite[Lemma 9.1.4]{JZ}. However, the condition  $K K^* - L L^* >0$ is in general not necessary, see \cite[Example 9.2.1]{JZ}.\\

It will now be shown that
the zero dynamics of systems of the form \eqref{eq:2.1d}--\eqref{eq:2.4d}  are again a port-Hamiltonian system, but with possibly a smaller state, that is, instead of $L^2(0,1;\mathbb C^n)$  the state space will be  $L^2(0,1;\mathbb C^k)$ with $0\le k\le n$.
First, it is shown that the problem of determining the zero dynamics for  \eqref{eq:2.1d}--\eqref{eq:2.4d} can be transformed into determining the zero dynamics for the finite-dimensional discrete-time system described by the matrices $A_d$, $B_d$, $C_d$ and $D_d$.
%
%
\begin{thm}
\label{T:4}
Let $z_0 \in L^2(0,1;{\mathbb C}^n)$. Then the following are equivalent.
\begin{enumerate}
\item There exists an input $u \in L^2_{loc}(0,\infty;{\mathbb C}^p)$ such that the output $y$ of \eqref{eq:2.1d}--\eqref{eq:2.4d}  with initial condition $z(\cdot, 0) = z_0$ is identically zero;
\item $z_0 \in L^2(0,1; V^*_d)$, where $V^*_d \subset {\mathbb C}^n$ is the largest output nulling subspace of the discrete-time system $\Sigma(A_d,B_d,C_d,D_d)$ with state space $\mathbb C^n$ given
by
\begin{equation}
\label{eq:23}
  w(n+1) = A_d w(n) +B_d u(n), \qquad y(n)= C_d w(n) + D_d u(n).
\end{equation}
\end{enumerate}
In particular, the largest output nulling subspace $V^*$ of \eqref{eq:2.1d}--\eqref{eq:2.4d} is given by $V^*=L^2(0,1; V^*_d)$.
\end{thm}
\begin{pf} 
The system \eqref{eq:2.1d}--\eqref{eq:2.4d} can be equivalently written as 
as \eqref{eq:19}--\eqref{eq:20}. In these equations the input, state and output were still spatially dependent. 
However,  the time axis has been split as  $[0,\infty) = \cup_{n \in {\mathbb N}} [n, (n+1)] . $
Thus condition 1.\ is equivalent to
\begin{enumerate}
\item[$1^\prime $] 
There exists a sequence $(u_d(n))_{n\in \mathbb 
N}\subseteq L^2(0,1; \mathbb C^m)$ and  a set $\Omega\subset (0,1)$ whose complement has measure zero such that for every $\zeta \in \Omega ,$ 
\begin{eqnarray}
   z_d(n+1)(\zeta) &=& A_d z_d(n)(\zeta) + B_d u_d(n)(\zeta), \label{z-discrete}\\
  (z_d(0))(\zeta) &=& z_0(1-\zeta). \nonumber\\
   0 &=& C_d z_d(n)(\zeta) + D_d u_d(n)(\zeta),  \nonumber 
\end{eqnarray}
\end{enumerate}
Clearly, condition $1^\prime$ implies that $z_0  (\zeta) \in V^*_d $ a.e., where $V^*_d$ denotes the largest output nulling subspace of the finite-dimensional system \eqref{eq:23}. Since trivially $z_0  \in L_2 (0,1; V^*_d),$  condition 2 follows.

The system $(A_d , B_d ,C_d, D_d)$ is a finite-dimensional discrete-time system. Let $V^*_d \subset \mathbb C^n$ indicate the largest output nulling subspace. Then there   exists a matrix ${\mathcal K}$ such that the output-nulling control is given by $u_d(n)= {\mathcal K} z_d(n)$, see \cite{wonham}.  Referring now to (\ref{z-discrete}),  if $z_0 \in L^2 (0,1; V^*_d) $  
then the output-nulling control $(u_d(n))_{n\in \mathbb 
N}$ for system (\ref{z-discrete}) satisfies $u_d(n)\in L^2 (0,1; \mathbb C^p).$ 
Condition   2 thus implies condition $1^\prime .$
\end{pf}
\medskip


For many partial differential equation systems, the largest output nulling subspace is not closed and the zero dynamics are not well-posed, \cite{MorrisRebarberIJC,Zwartbook}. However, for  systems of the form  \eqref{eq:2.1d}--\eqref{eq:2.4d} the  largest output nulling subspace is closed, and the zero dynamics are well-posed. 
The following theorem provides a characterization of the largest output nulling subspace of $\Sigma(A_d,B_d,C_d,D_d)$ and hence  of the zero dynamics for the original partial differential equation.  The proof can be found in \cite{IFAC2015}.
%
\begin{thm}
\label{T:5}
Define $E=-\left[\begin{smallmatrix} K_0\\ K_y \end{smallmatrix} \right]$, $F=\left[\begin{smallmatrix} L_0\\ L_y \end{smallmatrix} \right]$. 
The initial condition
$v_0$ lies in the largest output nulling subspace $V_d$ of $\Sigma(A_d,B_d,C_d,D_d)$ if and only if 
there exists a sequence $\{v_k\}_{k\geq 1} \subset {\mathbb C}^n$ such that
  \begin{equation}
    \label{eq:26}
    E v_{k+1} = F v_k, \qquad k\geq 0.
  \end{equation}
Furthermore, the largest output nulling subspace $V^*_d$ satisfies $V^*_d= \cap_{k \geq 0} V^k$, where $V^0={\mathbb C}^n$, $V^{k+1} = V^k \cap F^{-1}EV^k$.
\end{thm}

Thus in addition to the well-known $V^*$-algorithm for finite-dimensional systems, see \cite[p. 91]{Coron}, Theorem  \ref{T:5} provides 
an alternative  algorithm. It remains to show that the system restricted to the output nulling subspace is again port-Hamiltonian.
\begin{thm}
\label{T:6}
For the port-Hamiltonian system \eqref{eq:2.1d}--\eqref{eq:2.4d} the zero dynamics is well-posed, and the dynamics restricted to the largest output nulling subspace is a port-Hamiltonian system without inputs.
\end{thm}
\begin{pf}
By Theorem \ref{T:4},  the largest output nulling subspace $V^*$ of \eqref{eq:2.1d}--\eqref{eq:2.4d} is given by $V^*=L^2(0,1; V^*_d)$.
 If $V^*_d=\{0\}$, then there is nothing to prove, and so assume that $V^*_d$ is a non-trivial subspace of ${\mathbb C}^n$. It is well-known that there exists a matrix $F_d$ such that \cite{wonham}
\[
  (A_d + B_d F_d)V^*_d \subset V^*_d.
\]
Therefore,  using Theorem \ref{T:4} and \eqref{eq:19}--\eqref{eq:20}, it is easy to see that for the choice $u_d(n)(\zeta) := F_d z_d(n)(\zeta)$ the output $y_d(n)(\zeta)$ is  zero provided the initial condition $z_0$ lies in $L^2(0,1;V^*_d)$.
Using the definition of the $A_d$, $B_d$, $C_d$, $D_d$, $u_d$ and  $z_d$, it follows that for $z_0 \in L^2(0,1;V^*_d)$ there exists a function $f$ satisfying
\begin{eqnarray}
f(t) &=& z_0(1-t),\qquad t\in [0,1],\label{eq:18a+}\\
  \label{eq:18+}
   f( 1+ \lambda_0 t)  &=& - K^{-1} L  f( \lambda_0 t) -\lambda_0^{-1}K^{-1}  \begin{bmatrix} 0\\ F_d f(\lambda_0t) \end{bmatrix},\quad t\ge 0, \\
   0 &=& (\lambda_0K_y K^{-1} L-\lambda_0L_y) f(\lambda_0 t) + K_y K^{-1}  \begin{bmatrix} 0\\ F_d f(\lambda_0t) \end{bmatrix},\quad t\ge 0. \label{eq:18b+}
\end{eqnarray}
Equations \eqref{eq:18+}-\eqref{eq:18b+} can be equivalent written as
\begin{equation}
\label{eq:28}
  0 =  -\lambda_0 K_{ext} f(  1+\lambda_0 t )  -\lambda_0 L_{ext}   f(  \lambda_0  t ),
\end{equation}
with
\begin{equation}
\label{eq:29}
 K_{ext} = \begin{bmatrix} K\\ K_y \end{bmatrix}
\end{equation}
and some matrix $ L_{ext} $. 
Since $z_0 \in L^2(0,1;V^*_d)$,  for all $t$ and almost all $\zeta \in [0,1]$,  $f(\zeta + \lambda_0 t ) \in V^*_d$. Thus, $K_{ext}$ and $L_{ext}$ can be restricted to $V^*_d$ and  equation \eqref{eq:28} can equivalently be written with  matrices $K_{ext}|_{V^*_d}$ and $L_{ext}|_{V^*_d}$. Since $K$ is part of the  the matrix $K_{ext}$,  the matrix  $K_{ext}|_{V^*_d}$ has rank equal to the dimension of $V^*_d$. Let $P$ be the projection onto the range of $K_{ext}|_{V^*_d}$. 
This leads to \begin{equation}
\label{eq:30}
  0 =  -\lambda_0 PK_{ext}|_{V^*_d} f(  1+\lambda_0 t )  -\lambda_0 PL_{ext}|_{V^*_d}  f(  \lambda_0  t ).
\end{equation}
Define $K_{V^*_d}:=PK_{ext}|_{V^*_d}$ and $L_{V^*_d}:=PL_{ext}|_{V^*_d}.$ The above equation is the solution of the partial differential equation
\begin{eqnarray}
\label{eq:31}
  \frac{\partial }{\partial t}  z(\zeta,t)   &=&  -\lambda_0\frac{\partial}{\partial \zeta}  z(\zeta,t),\\
    \label{eq:32}
0 &=& -\lambda_0 K_{V^*_d}   z(0,t) -\lambda_0 L_{V^*_d}  z(1,t)
\end{eqnarray}
on the state space $L^2(0,1;V^*_d)$. Since $K_{V^*_d}$  is invertible, Theorem \ref{thm_JZ} implies that this  system is a  well-posed port-Hamiltonian system.
\end{pf}
\medskip

In the following section a present a second method to obtain the zero dynamics for systems with one dimensional input and output spaces is developed. The advantage of this method is that a transformation to a discrete system is not needed and  non-constant wave speed is possible.

\section{Zero dynamics of port-Hamiltonian systems with commensurate wave speed}

In this section the zero dynamics of systems of the form \eqref{eq:2.1d}-\eqref{eq:2.4d} with one dimensional input and output spaces and (possibly) non-constant wave speed are defined. The class of systems considered has the form
\begin{align}
\frac{\partial }{\partial t}    z(\zeta,t)  & = -\frac{\partial}{\partial \zeta} \left(\lambda_0(\zeta)  z(\zeta,t)\right) \label{eqnport1}\\
0   &= K_0 (\lambda_0(0)   z(0,t)) +  L_0  (\lambda_0(1)  z(1,t)) \label{eqnport2}\\
u(t)   &= K_u (\lambda_0(0)   z(0,t)) +  L_u  (\lambda_0(1)  z(1,t)) \label{eqnport3}\\
y(t)   &= K_y    z(0,t) +  L_y    z(1,t). \label{eqnport4}
\end{align}
Here $K_0, L_0\in \mathbb C^{(n-1)\times n}$, $K_u, K_y, L_u, L_y\in \mathbb C^{1\times n}$ and $\lambda_0\in L^\infty(0,1)$ satisfying $0<m\le \lambda_0(\zeta)\le M$ for almost every $\zeta\in(0,1)$ and constants $m,M>0$.
It will be assumed throughout this section that the port-Hamiltonian system \eqref{eqnport1}--\eqref{eqnport4} is a well posed linear system with state space $L^2(0,1;\mathbb C^n) $ or equivalently that the matrix $\left[\begin{smallmatrix} K_0\\K_u\end{smallmatrix}\right]$ is an invertible $n\times n$-matrix, see Theorem \ref{thm_JZ}. The corresponding generator $A$ of the $C_0$-semigroup of the homogeneous system is given by
\cite{JZ}
$$ Af= -(\lambda_0 f)', \qquad D(A)=\left\{\lambda_0 f\in H^1(0,1;\mathbb C^n)\mid \begin{bmatrix} 0 \\ 0 \end{bmatrix}  = \begin{bmatrix} K_{0}  \\ K_u  \end{bmatrix}   (\lambda_0 f)(0) + \begin{bmatrix} L_{0}  \\ L_u \end{bmatrix} (\lambda_0  f)(1) \right\}.$$

Denote by $G(s)$ the transfer function of the port-Hamiltonian system \eqref{eqnport1}--\eqref{eqnport4}. Since the port-Hamiltonian system is assumed to be well-posed, there exists a right half plane 
$$\mathbb C_\alpha:=\{s\in\mathbb C\mid {\rm Re}\, s>\alpha\}$$ 
such that $G: \mathbb C_\alpha\rightarrow \mathbb C$ is an analytic and bounded function.
Define 
\[ p:= \int_0^1 \lambda_0^{-1}(s) ds.\]
Moreover, using \cite[Theorem 12.2.1]{JZ} for $s\in\rho(A)$, where $\rho(A)$ denotes the resolvent set of $A$, and $u\in\mathbb C$ the number $G(s)u$ is (uniquely) determined by
\begin{align}
0 &= (K_0+L_0 e^{-s p}) v,\label{eqnzero1}\\
u &= (K_u+L_u e^{-s p}) v,\label{eqnzero2}\\
G(s)u&= (K_y+L_y e^{-s p}) v\label{eqnzero3}
\end{align}
for some $v \in \mathbb C^n.$
\begin{lem}\label{lemtransf}
There
  exists $\mu\in\mathbb R$ such that, for $s\in \mathbb C_\mu$,   $G(s)=0$ if and only if   the matrix
$\left[\begin{smallmatrix} K_0+L_0 e^{-sp} \\K_y+L_y e^{-sp}\end{smallmatrix}\right]$ is not  invertible.
\end{lem}
\begin{pf}
Since the matrix $\left[\begin{smallmatrix} K_0\\K_u\end{smallmatrix}\right]$ is invertible and $A$ generates a $C_0$-semigroup there is a $\mu\in\mathbb R$ such that $\rho(A)\subseteq \mathbb C_\mu$ and 
$$\left[\begin{array}{c} K_0+L_0 e^{-sp} \\K_u+L_u e^{-sp}\end{array}\right] = \left[\begin{array}{c} K_0\\ K_u \end{array}\right] + \left[\begin{array}{c} L_0 \\ L_u \end{array}\right]  e^{-sp}$$
is invertible for $s\in \mathbb C_\mu$. 

Assume now $G(s)= 0$ for some  $s\in \mathbb C_\mu$.  Then \eqref{eqnzero1}--\eqref{eqnzero3} imply that there exists $v \in\mathbb C^n$ such that 
\begin{align*}
0 &= (K_0+L_0 e^{-sp}) v,\\
1 &= (K_u+L_u e^{-sp}) v,\\
0&= (K_y+L_y e^{-sp}) v \, .
\end{align*}
Because $\left[\begin{smallmatrix} K_0+L_0 e^{-sp} \\K_u+L_u e^{-sp}\end{smallmatrix}\right]$ is invertible,  it yields $v \neq 0$. Thus $\left[\begin{smallmatrix} K_0+L_0 e^{-sp} \\K_y+L_y e^{-sp}\end{smallmatrix}\right]$ is not  invertible.

Conversely, assume that  for some $s\in\mathbb C_\mu$,  $\left[\begin{smallmatrix} K_0+L_0 e^{-sp} \\K_y+L_y e^{-sp}\end{smallmatrix}\right]$ is not  invertible. Then there exists a non-zero vector $v \in \mathbb C^n\backslash\{0\}$ such that
$$ \left[\begin{array}{c} 0\\0 \end{array}\right] =    \left[\begin{array}{c} K_0+L_0 e^{-sp} \\K_y+L_y e^{-sp}\end{array}\right] v.$$
Set $u:= (K_u+L_u e^{-sp})v$. Since $\left[\begin{smallmatrix} K_0+L_0 e^{-sp} \\K_u+L_u e^{-sp}\end{smallmatrix}\right]$ is invertible,  it follows that $u\not=0$. However, $G(s)u=0$ by \eqref{eqnzero1}--\eqref{eqnzero3}, which implies $G(s)=0$.
\end{pf}

\begin{thm}
\label{thm42}
 Suppose that $G(s)\not\equiv 0$. 
Then the zero dynamics of the port-Hamiltonian system \eqref{eqnport1}--\eqref{eqnport4}  are again a well-posed port-Hamiltonian system with wave speed $-\lambda_0 $ and  possibly a smaller state
space.  More precisely, there exists $k\in\{0,\cdots, n\}$ such that the zero dynamics is described by the port-Hamiltonian system 
\begin{align*}
\frac{\partial }{\partial t}    w(\zeta,t)  & = -\frac{\partial}{\partial \zeta} (\lambda_0 (\zeta)  w(\zeta,t)) \\
0   &= K_w    (\lambda_0 (0)w(0,t) )+  L_w    (\lambda_0 (1)w(1,t)).
\end{align*}
with state space $L^2(0,1;\mathbb C^k)$ and the $k\times k$-matrix $K_w$ is invertible.
\end{thm}
\begin{pf}
The zero dynamics are defined by  the  equations
\begin{align}
\frac{\partial }{\partial t}    z(\zeta,t)  & = -\frac{\partial}{\partial \zeta} (\lambda_0 (\zeta)  z(\zeta,t)) \label{eqnz1} \\ 
\begin{bmatrix} 0 \\ 0 \end{bmatrix}  &= \begin{bmatrix} K_{0}  \\ K_y  \end{bmatrix}   (\lambda_0 (0)z(0,t)) + \begin{bmatrix} L_{0}  \\ L_y \end{bmatrix}  (\lambda_0 (1)z(1,t) ) \, . \label{eqnz2}
\end{align}
Since there is one input and one output, and 
rank $ \left[\begin{smallmatrix} K_0\\K_u \end{smallmatrix}\right] = n$,  the  rank of the  matrix $\left[\begin{smallmatrix} K_0\\K_y\end{smallmatrix}\right]$ equals $n-1$ or $n.$

If rank $\left[\begin{smallmatrix} K_0\\ K_y \end{smallmatrix} \right] = n$, that is, this matrix is invertible, then the zero dynamics is well-defined on the whole state space $L^2(0,1;\mathbb C^n)$, see Theorem \ref{Tm:3}. Theorem \ref{thm_JZ} implies that the zero dynamics are well-posed on the state space $L^2(0,1;\mathbb C^n)$. Thus $k=n$ and the theorem is proved.

   Suppose next that rank $\left[\begin{smallmatrix} K_0\\ K_y \end{smallmatrix} \right]=n-1 .$ Then $K_y$ is a linear combination of the rows of $K_0$ and  there is an invertible transformation, a row reduction,  so that \eqref{eqnz2} is equivalent to
\begin{equation}\label{eqnz3}
\begin{bmatrix}
0 \\ 0 \end{bmatrix} =
 {\begin{bmatrix} K_{11} &K_{12}  \\ 0 &0 \end{bmatrix} } (\lambda_0 (0) z(0,t))   + \begin{bmatrix} L_{11} & L_{12}  \\  L_{21} & L_{22} \end{bmatrix}   (\lambda_0 (1) z(1,t)).
\end{equation}
Here $ K_{11},  L_{11} \in \mathbb C^{(n-1)\times (n-1)}$ and $   L_{22} \in \mathbb C$.
Since rank$[ K_{11}\,\, K_{12}] =n-1$,  column transformations lead to a representation where  the matrix $K_{11}$ is invertible. Assume now that this has been done.

Since  $K_{11}$ is invertible, and  $G $ is not equivalently zero,  Lemma \ref{lemtransf}, implies that 
   there exists $s_0\in\mathbb C$ such that  both $T_1:= K_{11}+L_{11} e^{-s_0p} $ and 
\begin{equation}  T:=\begin{bmatrix} T_1 & T_2  \\ T_3 &T_4 \end{bmatrix} := \begin{bmatrix} K_{11}+ L_{11} e^{-s_0p} & K_{12}+ L_{12} e^{-s_0p}  \\  L_{21} e^{-s_0p} & L_{22} e^{-s_0p} \end{bmatrix} \label{T}
\end{equation}
are invertible. 
Defining the Schur complement of $T$ with respect to $T_1$,
$$S= T_4-T_3 T_1^{-1} T_2 ,$$
$$ \begin{bmatrix} T_1 & T_2  \\ T_3 &T_4 \end{bmatrix} = \begin{bmatrix} I & 0  \\ T_3 T_1^{-1} &I \end{bmatrix}  \begin{bmatrix} T_1 & 0 \\0 &S \end{bmatrix}   \begin{bmatrix} I  & T_1^{-1}T_2  \\ 0 &I \end{bmatrix} .$$
 Since $T_1$ and $T$ are invertible, $S$ is invertible and 
$$  T^{-1}:=\begin{bmatrix} T_1^{-1}+T_1^{-1}T_2S^{-1}T_3T_1^{-1} & -T_1^{-1}T_2S^{-1}  \\ -S^{-1}T_3T_1^{-1}&S^{-1} \end{bmatrix}.$$

Now apply the state transformation 
$$ \tilde z = Tz.$$
The equations \eqref{eqnz3}  are equivalent to
\begin{align*}
\begin{bmatrix}
0 \\ 0 \end{bmatrix} =&
 {\begin{bmatrix} K_{11} &K_{12}  \\ 0 &0 \end{bmatrix} }  T^{-1}  (\lambda_0 (0)\tilde z(0,t))   + \begin{bmatrix} L_{11} & L_{12}  \\  L_{21} & L_{22} \end{bmatrix}    T^{-1} (\lambda_0 (1) \tilde z(1,t))\\
 =& \begin{bmatrix} K_{11}( T_1^{-1}+T_1^{-1}T_2S^{-1}T_3T_1^{-1})- K_{12}S^{-1}T_3T_1^{-1} &K_{12}S^{-1}-K_{11} T_1^{-1}T_2S^{-1}  \\ 0 &0 \end{bmatrix}    (\lambda_0 (0)\tilde z(0,t) )  \\
  &+ \begin{bmatrix} L_{11}( T_1^{-1}+T_1^{-1}T_2S^{-1}T_3T_1^{-1})- L_{12}S^{-1}T_3T_1^{-1} &L_{12}S^{-1}-L_{11} T_1^{-1}T_2S^{-1}  \\ 0 &e^{s_0p} \end{bmatrix}   (\lambda_0 (1)\tilde z(1,t)).
\end{align*}
Also the system of partial differential equations  \eqref{eqnz1} are equivalent to
 \begin{equation} \frac{\partial }{\partial t}     \tilde z(\zeta,t)   =-\frac{\partial}{\partial
\zeta}  (\lambda_0 (\zeta) \tilde z(\zeta,t)).
\label{eq-red}
\end{equation}
Thus, the transformed partial differential equation is identical to the original.
The general solution 
\[ \tilde z_n (\zeta, t) = \frac{c}{\lambda_0 (\zeta)} e^{\int_0^\zeta \lambda_0^{-1}(s) ds -t}\]
  and the boundary condition  $\tilde z_n (1,t) =0 $  imply that $\tilde z_n \equiv 0. $ 

Define 
\begin{equation}
w:= \left[\begin{smallmatrix} \tilde z_1 \\ \vdots \\ \tilde z_{n-1}\end{smallmatrix}\right],
\label{w}
\end{equation}
 and the  matrices
\begin{align*}
 K_w &:=K_{11}( T_1^{-1}+T_1^{-1}T_2S^{-1}T_3T_1^{-1})- K_{12}S^{-1}T_3T_1^{-1}\\
  L_w &:= L_{11}( T_1^{-1}+T_1^{-1}T_2S^{-1}T_3T_1^{-1})- L_{12}S^{-1}T_3T_1^{-1} \\
  K_{w12}&:=K_{12}S^{-1}-K_{11} T_1^{-1}T_2S^{-1}.
  \end{align*}
  Thus it yields
  $$ \begin{bmatrix} K_{11} &K_{12}  \\ 0 &0 \end{bmatrix}  T^{-1}  =  \begin{bmatrix} K_w & K_{w12} \\0 & 0 \end{bmatrix}.
    $$
Here  $K_{w12}$ is a $(n-1) \times 1 $-matrix and  rank$\,K_w \geq n-2. $
The zero dynamics is described by the reduced port-Hamiltonian system
\begin{align*}
\frac{\partial w}{\partial t} &=-   \frac{\partial }{\partial \zeta}(\lambda_0 w),\\
0 &= K_w  (\lambda_0 (0) w(0,t)) + L_w  (\lambda_0 (1)w(1,t)) . 
\end{align*}
The reduced system is well-posed on $L^2(0,1;\mathbb C^{n-1})$  if and only $K_w $ is invertible; that is, $K_w$ has rank $n-1$. If $K_w$ is invertible, then the theorem is proved.

Now suppose  that rank $\,K_w=n-2 .$
As in the first part, elementary row and column transformations can be used to put the boundary conditions for the reduced system into the form, again indicating the state variables by $w$,
$$
\begin{bmatrix}
0 \\ 0 \end{bmatrix} =
 {\begin{bmatrix} \tilde{K}_{11} &\tilde{K}_{12}  \\ 0 &0 \end{bmatrix} }  (\lambda_0 (0)w(0,t)  ) + \begin{bmatrix} \tilde{L}_{11} & \tilde{L}_{12}  \\  \tilde{L}_{21} & \tilde{L}_{22} \end{bmatrix}   (\lambda_0 (1) w(1,t)).
$$
where $\tilde{K}_{11}$  is invertible. 
Define
$$\tilde{T} (s) = K_w + L_w e^{-sp} . $$
In order to repeat the above procedure,  a complex number $s$ such that  $\tilde{T} $ and  $  \tilde{K}_{11}  +  \tilde{L}_{11} e^{-sp}$ are both invertible is needed.
Set  $s=s_0$. Define 
$$ X= T_1^{-1}+T_1^{-1}T_2S^{-1}T_3T_1^{-1}.$$ Recalling that $T_1=K_{11} + e^{-s_0 p} L_{11},$ $T_2=K_{12} + e^{-s_0 p} L_{12},  $
\begin{eqnarray*}
 K_w + L_w e^{-s_0 p } &=& K_{11} X-K_{12} S^{-1} T_3 T_1^{-1} +  e^{-s_0 p} L_{11} X - e^{-s_0 p} L_{12} S^{-1} T_3 T_1^{-1} \\
 &=&T_1 X -T_2S^{-1} T_3 T_1^{-1} \\
 &=&I+ T_2S^{-1}T_3T_1^{-1} -T_2S^{-1} T_3 T_1^{-1}\\
 &=& I.
\end{eqnarray*}
 Thus, with $s=s_0$, $ \tilde{T} (s) $ is invertible. 
Define
$$ f_w:\mathbb C_\alpha\rightarrow \mathbb C, \qquad f_w(s)= \det [ \tilde T (s) ] .$$
 and so  $f_w (s_0) = 1.$ Since $f_w$ is analytic, there is a sequence $s_n $, Re$s_n \to \infty$ with $f(s_n ) \neq 0 . $ Choose then $s_w$  so that  $  \tilde{K}_{11}  +  \tilde{L}_{11} e^{-s p }$ is invertible. 
Repeating the previous procedure leads to a port-Hamiltonian system with state-space $L^2(0,1;\mathbb C^{n-2}).$ Since each iteration leads to a state-space with fewer number of state variables,  this procedure is guaranteed to converge within $n$ steps. 
\end{pf}
\smallskip

Since the zero dynamics are a well-posed dynamical system, the following result is immediate.
\begin{corollary}
The invariant zeroes are contained in a left-hand-plane. 
\end{corollary}

One consequence of calculating the zero dynamics using the original port-Hamiltonian form is that it is easy to obtain the  input $u$ that zeroes the output.
Suppose only one state space reduction in Theorem \ref{thm42} is needed. The state space of the zero dynamics is $L^2(0,1;\mathbb C^{n-1})$. From (\ref{eqnport3}) and (\ref{T})
\begin{eqnarray*}
u(t) &=& K_u \lambda_0 (0) z(0,t) + L_u \lambda_0 (1) z(1,t)\\
&=& K_u \lambda_0 (0)  T^{-1} \tilde z(0,t) + L_u \lambda_0 (1) T^{-1} \tilde z(1,t). 
\end{eqnarray*}
In the zero dynamics, $\tilde z_n \equiv 0$. 
Defining $\tilde{K}_u $ to be the first $n-1$ columns of $K_u \lambda_0 (0)  T^{-1}$ and defining $\tilde L_u$ similarly,
the zeroing input is 
$$ u(t) =\tilde K_u w (0,t) + \tilde L_u w(1,t)$$
where $w$ is defined in (\ref{w}).
For the situation where more than one state space reduction is needed, the calculation is similar, except that a transformation matrix $T$ is needed for each reduction.
\section{Computation}
Theorem 4.2 leads to a characterization of the zero dynamics as a port-Hamiltonian system of smaller dimension. Moreover, the proof is constructive and can be used in an algorithm  to calculate the
zero dynamics  using standard linear algebra algorithms, see the box on the following page. Zero dynamics can be calculated exactly  for large system order; that is those with a large number of nodes.  Furthermore, Theorem \ref{thm32} can be used to check stability. 

Several examples are now presented to illustrate
the calculation of zero dynamics.
\begin{eg}
Consider the system from Example \ref{E:3.1}, written in the equal wave speed form.
For zero dynamics,
\begin{equation}
\begin{bmatrix}
0 \\ 0 \\ 0 \end{bmatrix} = \underbrace{\begin{bmatrix} 1& 0 &1 \\ 0 & 1 & 0 \\ 0 & 0 &0 \end{bmatrix} }_{\Large K}  z(0,t) + \underbrace{\begin{bmatrix} 1 & 0 & 0 \\ 0 & 0 & -1 \\ 1 & 1 & 0 \end{bmatrix}}_{\Large L} z(1,t) .
\label{eg-eg1}
\end{equation}
The rank of $K=2$ and so the zero dynamics are defined on a smaller state space than the original. Applying one iteration of  the algorithm  yields (with $s_0=0$)
$$TP=\begin{bmatrix} 2 & 0 &1 \\ 0 & 1 & -1 \\1 & 1 & 0 \end{bmatrix} , \quad K_w =\begin{bmatrix} 0 & -1 \\ -1 & -1\end{bmatrix} , \quad L_w = \begin{bmatrix} 1 & 1 \\ 1 & 2 \end{bmatrix} .$$
The last row of the transformation matrix  $TP$ indicates that for zero dynamics
$$z_1 + z_{2a} \equiv 0 $$ and  the first two rows define the remaining state variables:
$$\tilde{z}_1 = 2 z_1 + z_{2b}, \quad \tilde{z}_2 = z_{2a}-z_{2b}. $$
(Because $z_1 + z_{2a} \equiv 0$, $\tilde z_1 $ can also be written $-2 z_2a + z_{2b} .$)
The rank $K_w=2. $ The zero dynamics are well-defined on $L^2(0,1;\mathbb C^{2})$ with governing differential equation (\ref{eq-red}) and boundary conditions
$$\begin{bmatrix} 0 & -1 \\ -1 & -1\end{bmatrix}  \tilde{z} (0,t) +  \begin{bmatrix} 1 & 1 \\ 1 & 2 \end{bmatrix}  \tilde{z} (1,t) .$$

From the definition of the control system in Example \ref{E:3.1},
$$u(t) = z_{2b} (0,t)$$
is the zeroing input.

This example could be done by hand. The definition of the zero dynamics (\ref{eg-eg1}) implies that $z_1+z_{2a} =0 .$ Substitution into the equations yields the zero dynamics.
\end{eg}

\noindent
\fbox{
\begin{minipage}{\textwidth}
\vspace{1.5ex}
\begin{center} {\bf  \underline{ Algorithm:   Calculation of Zero Dynamics} }\\[2ex]
\end{center}

The data are: wave speed $p=\int_0^1 \frac{1}{\lambda_0 (\xi ) } d \xi $, boundary condition matrices $K_0$, $L_0$, and output  matrices $K_y, $ $L_y$. 
The dimension of the system is $n$, the number of columns in $K_0.$
Define
$$K=\begin{bmatrix} K_0\\K_y\end{bmatrix}.$$
If $K$ is invertible
the zero dynamics are well-defined with $n$ state variables. Otherwise do the following calculations.
\begin{enumerate}
\item
Perform LU-decomposition of $K$: $P_{\ell u} K=M_\ell M_u $ where $M_\ell$ is lower triangular, $M_u$ is upper triangular and $P_{\ell u}$ is a permutation matrix.
\item
If necessary permute last column  of $M_u$ with earlier column so that rank of top left $n-1$ block is $n-1$; call the permutation matrix $P.$

Partition $M_u P$ and $ M_\ell^{-1} P_{\ell u} L P$ similarly as
$$M_u = \begin{bmatrix} K_{11} & K_{12} \\ \begin{bmatrix} 0  \ldots 0 \end{bmatrix} & 0 \end{bmatrix} , \quad M_\ell^{-1} P_{\ell u}  L = \begin{bmatrix} L_{11} & L_{12} \\ L_{21} & L_{22} \end{bmatrix} \, . $$

\item
Define the matrices
\begin{eqnarray*}
T_{1}= K_{11}+L_{11} e^{-p s_0} , \quad 
T=\begin{bmatrix} T_{1} & K_{12}+L_{12} e^{-p s_0} \\ L_{21} e^{-p s_0} & L_{22} e^{-p s_0} \end{bmatrix} 
\end{eqnarray*}
for  $s_0$ so that both matrices are invertible. (The existence of such an $s_0$ is guaranteed if the transfer function is not identically zero. A simple way find a suitable $s_0$ is to start with $s_0=0$ and then increase by an arbitrary amount until both matrices are invertible. )
\item
Decompose $T^{-1}$ using the same decomposition as for $K_u$ and construct the inverse of $T$ using the Schur complement. Letting $X$ be the solution of 
$ T_3 = X T_1, $ define
$$S_i=(T_4- X T_2 )^{-1}.$$
 (Note $S$ is a scalar.) Only the 2 left blocks of $T^{-1}$ are needed: 
$$(T^{-1})_{11} =T_1^{-1} (I+  T_2 S_k X) , \quad (T^{-1})_{21} = -S_i  X. $$
\item
The boundary matrices for the reduced system are
$$K_w=K_{11} (T^{-1})_{11} +K_{12} (T^{-1})_{21} , \quad  L_w=L_{11} (T^{-1})_{11} +L_{12} (T^{-1})_{21}      . $$
\item
The new variables are $\tilde{z}_1 \ldots \tilde{z}_{n-1}$ where $\tilde{z}=T P z $, the differential equation is
$$\frac{\partial }{\partial t}     \tilde z(\zeta,t)   =-\frac{\partial}{\partial
\zeta}  (\lambda_0 (\xi)  \tilde z(\zeta,t)) $$
and the boundary conditions are
$$K_w \lambda_0 (0) \tilde{z } (0,t) + L_w \lambda_0 (1) \tilde{z} (1,t) .$$
If rank $K_w=n-1$, the algorithm is complete. 
If not, return to the first step with $K=K_w$, $L=L_w$ and repeat the process.
 \end{enumerate}
\vspace{.5ex}
\end{minipage}
}\smallskip
\newpage

%
%
\begin{eg}
$$ \frac{ \partial x_i}{\partial t} = - \frac{ \partial x_i}{\partial \zeta }, \hspace{2em} i=1,2,3 . $$
with 
\begin{eqnarray}
\label{eg2-u} 
\begin{bmatrix}0 \\ 0 \\ u (t) \end{bmatrix} &=& \underbrace{\begin{bmatrix} 0 & 0 &-1 \\ 0 & -1 & 0 \\ 0 & -1 & 0 \end{bmatrix}}_{\Large K} x(0,t) + \begin{bmatrix} 1 & 0 & 0 \\ 0 & 1& 0 \\ 0 & 0  &1 \end{bmatrix} x(1,t) \\
y(t) & = & \begin{bmatrix} 0 & 0 & 0 \end{bmatrix} x(0,t) + \begin{bmatrix} 1 & 0 & 0 \end{bmatrix} x(1,t) . \nonumber
\end{eqnarray}
The rank of $K$ in (\ref{eg2-u}) is 3 and so the system is well-posed. The transfer function is not identically zero. 

Zero dynamics require
\begin{eqnarray}
\begin{bmatrix}
0 \\ 0 \\ 0 \end{bmatrix} &=& \begin{bmatrix} 0 & 0 &-1 \\ -1 & 0 & 0 \\ 0 & 0 &0 \end{bmatrix} x(0,t) + \begin{bmatrix} 1 & 0 & 0 \\ 0 & 1 & 0 \\ 1 & 0 & 0 \end{bmatrix} x(1,t) .
\label{eg2:1}
\end{eqnarray}
%
Applying   the algorithm  yields (with $s_0=0$)
$$TP=\begin{bmatrix} -1& 1 &0 \\ 1 & 0 & -1 \\1 & 0 & 0 \end{bmatrix} , \quad K_w =\begin{bmatrix} 0& 0 \\ 0 & 1\end{bmatrix} , \quad L_w = \begin{bmatrix} 1& 0 \\ 0 & 0 \end{bmatrix} .$$
The third row of $TP $ implies that $z_1\equiv 0.$ The reduced states are 
$$\tilde z_2 = -z_1 + z_2 = z_2 , \quad \tilde z_3 = z_1 -z_3 =-z_3 .$$

Since  $K_w$ does not have full rank. the algorithm needs to be repeated; but with $K_w$, $L_w$ as the boundary matrices. This yields
$$(TP)_2=\begin{bmatrix} 0 & 1 \\ 1 & 0 \end{bmatrix} , \quad (K_w)_2 = \begin{bmatrix} 1 \end{bmatrix} , \quad (L_w)_2 = \begin{bmatrix} 0 \end{bmatrix} . $$
Thus $\tilde z_2 = z_2 \equiv 0$ and $\tilde{z}_3  (0)= -z_3 (0) = 0.$

This example is also simple enough to do by hand. The original equations  (\ref{eg2:1}) are already row-reduced, and imply $x_1 \equiv 0 $. The reduced system must have  $x_2 \equiv 0. $  

Either calculation  leads to one non-zero equation, for $x_3 $ with the boundary condition
$$x_3 (0,t) = 0. $$
The system equations \eqref{eg2-u} imply that in order to achieve this, $u (t) = x_3 (1,t).$
%
\end{eg}

\begin{eg}
Consider a larger system with $n=10.$ Suppose the wave speed $\lambda_0$ is such that $-\int_0^1 \lambda_0 (\xi ) d \xi = -1.$
The entries in the boundary matrices are zero, except that 
\begin{alignat*}{6}
K_0(1,2) &=  1,& \; &K_0(1,9)=-3 ,& \; &K_0(2,3)=1, \, &K_0(2,2)=-1, \;  &K_0(3,6)=1 ,  \; K_0 (3,10)=2,  \\
K_0(4,1) &= -5,&\; &K_0(4,6)=2,&\;  &K_0(5,10) =  6 , \; &K_0(5,9)=-4, \; &K_0(6,8)=4, \; K_0(6,1)=-2, \\
K_0 (7,6)&= 1,& \; &K_0(7,7)=3,& \; &K_0(8,3)=-2, \; &K_0(8,8)= 1, \; &K_0(8,5)=-5, \\ 
K_0(9,1)&=1,& \; &K_0(9,6)=5,& \; &K_0(9,9)=-1 ;&&\\
K_u(1,4)&= 1 ;& &&&&&\\
L_y(1,2) & = 1,& \; &L_y (1,4)=-2. &&&&
\end{alignat*}
Since
$${\rm rank } \begin{bmatrix} K_0\\ K_u \end{bmatrix} = 10 $$
this system is well-posed.  Also, the transfer function $G$ is not identically $0$; in particular $G(0) \neq 0 .$
Applying the algorithm with  $s_0 = 0$ yields
\begin{eqnarray*}
TP &=& \begin{bmatrix}
   -5  & 0    &      0  &        0  &        0    & 2 & 0   &       0  &        0     &     0\\
         0    &-1  &   1        &  0   &       0  &        0    &      0   &       0    &      0   &       0\\
         0         & 0   & -2  &        0  &  -5   &       0   &       0   &  1  &       0         & 0\\
         0      &    0   &       0  &        0  &  -2.5  &        0     &     0  &   0.5  &   -3          &0\\
         0        &  0   &       0      &    0    &      0   &       0    &      0     &     0 &   -4 &    6\\
         0       &   0      &     0  &        0   &       0  &   5.4  &       0     &     0   & -1  &         0\\
         0     &     0    &      0   &       0    &      0   &       0  &   3  &         0    & 0.1852      &    0\\
         0     &     0    &      0    &      0    &      0    &      0   &       0   &  4  &  -0.1481      &    0\\
         0     &     0     &     0     &     0    &      0     &     0    &      0     &     0 &   0.1852   &       0\\
         0    & 1  &        0  & -2 &         0    &      0    &      0    &      0   &       0    &      0 
         \end{bmatrix} 
         , \quad
K_w = I_9, \quad L_w = 0_{9 \times 9 }.
   \end{eqnarray*}
   For zero dynamics,
$z_2 -  2 z_4 \equiv 0$ and the zeroing input is
$$u(t) =K_u TP z(0,t) = -2.5 z_5 (0,t) +0.5 z_8 (0,t)-3 z_9 (0,t).$$
\end{eg}

\section{Conclusions}

In this paper, zero dynamics were formally defined for port-Hamiltonian systems. If the feedthrough operator is invertible, then the zero dynamics are again a port-Hamiltonian system of the same order. In general, however, the feedthrough operator is not invertible. For many  infinite-dimensional systems, where the feedthrough is not invertible, the zero dynamics  are not well-defined. It has been  shown in this paper that provided the system can be rewritten as a network of waves with the same speed, the zero dynamics are always well-defined, and are a port-Hamiltonian system. Furthermore, a numerical method to construct the zero dynamics using the original partial differential equation has been described. Finite-dimensional approximations, which can be inaccurate in calculation of zeros, are not needed. The approach applies to systems with commensurate but non-equal wave speeds, and this generalization will be explored in future work. The extension to multi-input multi-output systems also needs to be established.


\end{document}